\documentclass{gtart}
\usepackage{amsfonts,amssymb}
\usepackage{amsmath,amscd}

\input gtmonout
\volumenumber{2}
\volumeyear{1999}
\volumename{Proceedings of the Kirbyfest}
\pagenumbers{125}{133}
\papernumber{7}
\received{30 December 1998}\revised{3 July 1999}
\published{17 November 1999}

\newtheorem{thm}{Theorem}[section]

\newtheorem{prop}[thm]{Proposition}

\def\mc{\mathcal}
\def\B{\mathbb B}

\def\L{\mathbb L}
\def\N{\mathbb N}

\def\R{\mathbb R}
\def\S{\mathbb S}

\def\*{\times}
\def\al{\alpha}
\def\be{\beta}

\def\LA{\Lambda}
\def\om{\omega}
\def\la{\langle}
\def\ra{\rangle}
\def\nix{\varnothing}

\begin{document}
\title{Manifolds at and beyond the limit of metrisability} 
\authors{David Gauld}

\address{Department of Mathematics, The University of Auckland
\\Private Bag 92019, Auckland, New Zealand}
\email{gauld@math.auckland.ac.nz}

\begin{abstract}
In this paper we give a brief introduction to criteria for
metrisability of a manifold and to some aspects of non-metrisable
manifolds. Bias towards work currently being done by the author and
his colleagues at the University of Auckland will be very evident.

{\em Wishing Robion Kirby a happy sixtieth with many more to follow.} 

\end{abstract}
\asciiabstract{In this paper we give a brief introduction to criteria
for metrisability of a manifold and to some aspects of non-metrisable
manifolds. Bias towards work currently being done by the author and
his colleagues at the University of Auckland will be very evident.}

\maketitle

\section{Introduction}

By a \emph{manifold} we mean a connected, Hausdorff topological space in which each point has a neighbourhood homeomorphic to some euclidean space $\R^n$. Thus we exclude boundaries. Although the manifolds considered here tend, in a sense, to be large, it is straightforward to show that any finite set of points in a manifold lies in an open subset homeomorphic to $\R^n$; in particular each finite set lies in an arc. 

Set Theory plays a major role in the study of large manifolds; indeed, it is often the case that the answer to a particular problem depends on the Set Theory involved. On the other hand Algebraic Topology seems to be of less use, at least so far, with \cite[Section 5]{N} and \cite{DG} being among the few examples. Perhaps it is not surprising that Set Theory is of major use while Algebraic Topology is of minor use given that the former is of greater significance when the size of the sets involved grows, while the latter is most effective when we are dealing with compact sets. Thus a large part of this study involves the determination of relationships between various topological properties when confined to manifolds. Quite a lot of effort has gone into determining whether a particular combination of properties is equivalent to metrisability. 

The theory of non-metrisable manifolds is much less developed than that of compact manifolds, but again this is not surprising for, as noted in Section 3, there are many more, $2^{\aleph_1}$, of the former than the latter. The nearest there is to a classification of non-metrisable surfaces is the Bagpipe Theorem of Nyikos, \cite{N}, which states that any
$\om$--bounded surface is the union of a compact subsurface together with a finite number of mutually disjoint ``long pipes.'' A space is \emph{$\om$--bounded} if every countable subset has compact closure, and a \emph{long pipe} is a
manifold-with-boundary which is the union of an increasing $\om_1$--sequence of open subsets each of which is homeomorphic to $\S^1\* [0,1)$. The paper \cite{N} gives an excellent introduction to the theory of non-metrisable manifolds, while \cite{N1} provides a more recent view.

Although Set Theory had been used in the study of manifolds earlier, the solution of the following important question (\cite{A} and \cite{W}) marks a watershed in its use: \begin{itemize}
\item Must every perfectly normal manifold be metrisable? \end{itemize}
A topological space $X$ is \emph{perfectly normal} provided that for each pair of disjoint closed subsets $A,B\subset X$ there is a continuous function $f\co X\to [0,1]$ such that $f^{-1}(0)=A$ and $f^{-1}(1)=B$. In the 1970s it was found that the answer to this question depends on the Set Theory being used. More precisely, Rudin and Zenor in \cite{RZ} constructed a counterexample in Set Theories satisfying ZFC (=Zermelo--Fraenkel plus the Axiom of Choice) and the Continuum Hypothesis whereas Rudin in \cite{R} proved that in ZFC Set Theories satisfying Martin's Axiom and the negation of the Continuum Hypothesis every perfectly normal manifold is metrisable. 

One might identify three levels at which Set Theory is applied to the study of manifolds. At the basic level is the application of properties which are valid in all Set Theories (say those satisfying ZFC): these include the Well-Ordering Principle (frequently applied in the form of transfinite induction) and the Pressing Down Lemma. At the next level is the application of standard axioms which are known to be independent of ZFC: these include the Continuum Hypothesis and Martin's Axiom. At the third level is the application of such techniques as forcing to construct Set Theories satisfying certain properties (of course that is how use of the Continuum Hypothesis etc is justified, but someone else has already done the hard Set Theory). 

Denote by $\om_1$ the collection of countable ordinals: use the Well-Ordering Principle to well-order an uncountable set and then declare $\om_1$ to be the subset consisting of those having only countably many predecessors. We can topologise $\om_1$ by using the order topology. Closed unbounded subsets in this topology play an important role. Any subset of $\om_1$ which meets every closed unbounded set is called a \emph{stationary} set; for example the set of all limit points of $\om_1$, indeed any closed unbounded subset, is stationary. If $S\subset\om_1$ then a function $f\co S\to\om_1$ is called \emph{regressive} provided that $f(\al )<\al$ for each $\al\in S$. The proof of the following proposition may be found in many books on Set Theory.

\begin{prop}[Pressing Down Lemma] Let $f\co S\to\om_1$ be a regressive function, where $S$ is stationary. Then there are $\be\in\om_1$ and a stationary set $T\subset S$ such that $f(T)=\{\be\}$. \end{prop}

\noindent {\bf Continuum Hypothesis, CH}\qua (Cantor)\qua {\sl Any subset of $\R$ either has the same cardinality as $\R$ or is countable.}
\bigskip

\noindent {\bf Martin's Axiom, MA}\qua {\sl In every compact, ccc, Hausdorff space the intersection of fewer than $2^{\aleph_0}$ dense open sets is dense.}
\bigskip

A space has the \emph{countable chain condition} (abbreviated ccc) provided that every pairwise disjoint family of open sets is countable. It is well-known that CH is independent of ZFC. It is also the case that MA is independent of ZFC but from the Baire Category Theorem it is immediate that CH$\Rightarrow$MA. Thus we might expect three possible kinds of ZFC Set Theories in which combinations of CH and MA appear: ZFC+CH, ZFC+MA+$\neg$CH and ZFC$+\neg$MA. All three possibilities do occur. Usually in applications other equivalent versions of MA are used. 

\section{Conditions Related to Metrisability} \cite{G} contains about 50 conditions equivalent to metrisability for a manifold. This does not include conditions which are equivalent in a general topological space. Nor does it include conditions which are equivalent only in some Set Theories. Some of the following conditions will be quite familiar, others may be less so. \begin{thm}\label{m}
The following conditions are equivalent for a manifold $M$: \begin{enumerate}
\item $M$ is metrisable;
\item $M$ is paracompact;
\item $M$ is nearly metaLindel\"of;
\item $M$ is Lindel\"of;
\item $M$ is second countable;
\item $M$ is finitistic;
\item $M\* M$ is perfectly normal;
\item $M$ may be properly embedded in some euclidean space; \item there is a cover $\mc U$ of $M$ such that for each $x\in M$ the set st$(x,\mc U)$ is open and metrisable;
\item the tangent microbundle on $M$ contains a fibre bundle. \end{enumerate}
\end{thm}

Proof of the equivalence of these conditions will not be given here, refer instead to \cite{G} where there is a discussion and references. However, to illustrate the first level of application of Set Theory here is a sketch of the proof of 9$\Rightarrow$4; full details appear in \cite{GG}. 

Let $\mc U$ be as in 9 and suppose $M$ has dimension $m$. For each $\al\in\om_1$ define inductively Lindel\"of, connected, open subsets $V_\al\subset M$. Let $V_0$ be any open subset of $M$ homeomorphic to $\R^m$. If $\lambda$ is any non-zero limit ordinal let $V_\lambda =\cup_{\al <\lambda}V_\al$. If $\al$ is any ordinal and $V_\al$ is defined then $V_\al$ is separable; let $D_\al$ be a countable dense subset. Then $\overline V_\al\subset\cup_{d\in D_\al}st(d,\mc U)$, hence we may let $V_{\al +1}$ be the component of $\cup_{d\in D_\al}st(d,\mc U)$ containing $\overline V_\al$. Note that $\overline V_\be\subset V_\al$ whenever $\be <\al$.

Connectedness of $M$ implies that $M = \bigcup _{\alpha \in \omega _1} V_\alpha$. (We need to show
that $\bigcup _{\alpha \in \omega _1} V_\alpha$ is both open and closed. If $x$ is in the closure of
this set then the $n$th member of a countable neighbourhood base at $x$ meets some $V_{\beta_n}$ and
hence every member of this neighbourhood base meets $V_\beta$ where $\beta$ is the supremum of the
ordinals $\beta_n$. Thus $x\in\overline{V_\beta}$ so $\bigcup _{\alpha \in \omega _1} V_\alpha$ is
closed.) If for some limit ordinal $\lambda$ we have $\overline V_\lambda - V_\lambda =\nix$ then
$V_\lambda$ is open and closed, so is all of $M$ (by connectedness) and hence $M$ is Lindel\"of.

Suppose instead that for each limit ordinal $\lambda$ we have $\overline V_\lambda - V_\lambda \not= \nix$. Then $M$ cannot be Lindel\"of but we will obtain a contradiction. Choose $x_\lambda \in \overline V_\lambda - V_\lambda$, let $\LA$ denote the set of limit ordinals excluding 0 and define the regressive function $f\co \LA\to\omega_1$ by
\[ f(\lambda )=\mbox{min}\{\alpha\in \omega_1: st(x_\lambda ,{\mathcal U})\cap V_\alpha \not= \nix \}.\]
Then by the Pressing Down Lemma we may choose $\alpha\in\omega_1$ such that $A=f^{-1}(\alpha )$ is stationary.

For each $\lambda \in A$ choose $d_\lambda \in st(x_\lambda ,{\mathcal U})\cap D_\alpha$ and let $cst_\lambda (d_\lambda ,\mathcal U)$ denote the (open!) component of $st(d_\lambda ,\mathcal U)$ containing $x_\lambda$. Define $g\co A\to\omega_1$ by \[ g(\lambda )=\mbox{min}\{\beta \in \omega_1: cst_\lambda (d_\lambda ,{\mathcal U})\cap V_\beta \not= \nix \}.\] Again by the Pressing Down Lemma there is $\beta \in \omega_1$ so that $B=g^{-1}(\beta )$ is stationary.

Because $cst_\lambda (d_\lambda ,\mathcal U) \cap D_\beta\not= \nix$ for each $\lambda \in B$, and $D_\alpha$ and $D_\beta$ are countable whereas $B$ is uncountable, we may choose $d\in D_\alpha$ and $d'\in D_\beta$ such that
\[ C=\{\lambda\in B\ :\ d_\lambda =d \mbox{ and } d'\in cst_\lambda (d_\lambda ,\mathcal U)\}\]
is uncountable. One can show that $N=\bigcup _{\lambda \in C}cst_\lambda (d,\mathcal U)$ is connected (in fact
$N=cst_\lambda (d,\mathcal U)$ for each $\lambda \in C$) and hence $N$ is a metrisable manifold, so is hereditarily
Lindel\"{o}f.

Let $X=\{ x_\lambda \in M\ :\ \lambda \in C\}$. On the one hand, because $X\subset N$ and $N$ is hereditarily Lindel\"{o}f, $X$ is Lindel\"{o}f. On the other hand, $\{V_\lambda : \lambda\in\Lambda\}$ forms an open cover of $X$ with no countable subcover, and hence $X$ cannot be Lindel\"{o}f, a contradiction.\hfill$\bullet$ 

While not all terms introduced in Theorem \ref{m} will be defined here, it is worth noting that the term \emph{finitistic} was introduced by Swan, \cite{S}, in his study of the theory of fixed point sets. In \cite{M}, Milnor introduced the notion of microbundle to provide a means of adapting techniques developed for smooth manifolds to the topological context. Much of the impetus was lost when Kister in \cite{K} showed that many microbundles are equivalent to fibre bundles. Condition 10 of Theorem \ref{m} shows that Kister's result goes no further than metrisable manifolds when attention is fixed on the tangent microbundle, which was one of Milnor's prime concerns. 

It is interesting to note that the rather strong condition of second countability is equivalent to metrisability. Consequently every metrisable manifold is separable. Perhaps this is not surprising in view of the fact that there is a kind of duality between separability and Lindel\"ofness. However it is not hard to construct separable manifolds which are not metrisable, for example the Pr\"ufer manifold of the next section. Furthermore we have the same set theoretic dilemma as for Alexandroff's question: in ZFC$+$CH one can construct (see \cite{G2}) a non-metrisable manifold such that all of its finite powers are hereditarily separable; the manifold cannot be Lindel\"of. On the other hand Kunen has shown in \cite{Ku} that in ZFC$+$MA$+\neg$CH every space all of whose finite powers are hereditarily separable is Lindel\"of, hence metrisable if a manifold.

\section{Examples of Non-metrisable Manifolds} There are only 4 distinct 1--manifolds: the circle, $\S^1$, is the only one which is compact; the real line, $\R$, is the only one which is metrisable but not compact; then there are two non-metrisable 1--manifolds, the long line, $\L$, and the open long ray, $\L_+$. The latter two are distinguished by the fact that $\L$ (like $\S^1$ but unlike $\R$ and $\L_+$) is
$\om$--bounded. The \emph{long line}, $\L$, \cite{C}, is constructed as follows. Topologise $\om_1\* [0,1)$ using the lexicographic order, which is defined by $(\al ,s)<(\be ,t)$ provided either $\al <\be$ or else $\al =\be$ and $s<t$. Intervals of the form $(a,b)=\{ c\in\om_1\* [0,1)\ /\ a<c<b\}$ together with those of the form $[0,b)$ (usual meaning for half-open intervals; we have abbreviated the first element $(0,0)$ to 0) define a basis for the topology of $\om_1\* [0,1)$, which thereby becomes the closed long ray. Then $\L$ is obtained by joining together two copies of $\om_1\* [0,1)$ at 0, and $\L_+$ is a single copy with 0 removed.

When it comes to 2--manifolds, the situation is almost helpless. On page 669 of \cite{N} it is shown that there are at least $2^{\aleph_1}$ distinct $\om$--bounded, simply connected 2--manifolds, $2^{\aleph_1}$ being the cardinality of the power set of the first uncountable ordinal $\aleph_1$. It is quite a challenge to find enough topological invariants to identify this many distinct surfaces so, as one might expect, this part of the paper makes interesting reading. From \cite{H} it follows that there are exactly $2^{\aleph_1}$ distinct manifolds; contrast this with the situation for compact manifolds where, by \cite{CK}, there are only countably many. 

One of the simplest non-metrisable surfaces is the Pr\"ufer manifold, one version of which is the following. Let
\[ A = \{(x,y) : x,y\in \R \mbox{ and } x\not=0 \} \mbox{ and for each } y\in \R \mbox{, let } B_y =\{(0,y)\}\*\R .\] Write $\langle y,z \rangle$ for a typical member $((0,y),z)$ of $B_y$. Set \mbox{$M = A\cup
(\cup_{y\in\R} B_y)$ .} To make it into a manifold, $M$ is topologised by using the usual topology on $A$ and replacing each point, $(0,y)$, of the missing $y$--axis by the copy, $B_y$, of the real line. More precisely, suppose that $\langle y,z \rangle \in B_y$. Declare
\[ \{ \{ \langle y,\zeta \rangle \in B_y : |\zeta -z|<\frac{1}{n} \} \cup T(y,z,n) : n\in \N \}\]
to be a neighbourhood basis at $\langle y,z \rangle$, where \[ T(y,z,n) = \{ (\xi ,\eta )\in \R^2 :0<|\xi |<\frac{1}{n} \mbox{ and } z-\frac{1}{n}<\frac{\eta -y}{|\xi
|}<z+\frac{1}{n} \} .\]
That $M$ is a 2--manifold follows from the fact that the two triangles forming $T(y,z,n)$ may be opened out to make room for the interval $\{ y\}\* (z-1/n,z+1/n)$ to give a region homeomorphic to a square. Then $M$ is not metrisable because it is not Lindel\"of, containing an uncountable closed discrete subset. 

In \cite{GGKM} and \cite{Mo} two modifications are given to this construction, both giving rise to non-metrisable surfaces which satisfy quite strong conditions. The authors begin with two disjoint subsets $C, D\subset\R$ and for each $x\in D$ a sequence $\la x_n\ra$ in $C$ converging to $x$. They let $X(C,D)$ be the union $C\cup D$ but with points of $C$ isolated and neighbourhoods of points of $D$ containing tails of the corresponding sequences. The essential features of $X(C,D)$ are then preserved in a surface constructed along the lines of the Pr\"ufer manifold. The idea is essentially to take standard Pr\"ufer neighbourhoods of points of $A$ and $B_y$ whenever $y\in C$ while neighbourhoods of points of $B_y$ when $y\in D$ also include connected tails of Pr\"ufer neighbourhoods. In one case $C$ is a Bernstein subset of $\R$ (ie, every uncountable closed subset of $\R$ meets both $B$ and $\R -B$); such a set may be constructed using the Axiom of Choice so the manifold exists in ZFC. The resulting manifold possesses nice metric-type properties but is not perfect, hence not metrisable. In the other case forcing is used to construct the set $C$ satisfying the extra condition of countable metacompactness, which then carries over to the manifold.

Reference has already been made to two surfaces which have been constructed by use of the Continuum Hypothesis. Actually that found in \cite{G2} is a refinement of that in \cite{RZ}. The basic idea is to start with the \emph{open} unit disc $\B^2$ together with a half-open interval allocated to each point of the boundary $\S^1$. Using CH we may identify the points of $\S^1$ with $\om_1$, so the intervals are denoted $I_\al$ for $\al\in\om_1$. The manifold has $M=\B^2\cup (\cup_{\al\in\om_1}I_\al )$ as the underlying set and the topology is obtained from the usual topology on $\B^2$ by inserting the interval $I_\al$ into the boundary of $\B^2$ at the point $\al$. Let $M_\al =\B^2\cup (\cup_{\be <\al}I_\be )$. The insertion is done by induction on $\al\in\om_1$. Before beginning, the countable subsets of $M$ are also indexed by $\om_1$: under CH there are exactly as many countable subsets of $M$ as there are points of $\om_1$ so we may denote the $\al$th countable set by $S_\al$. In much the same way as the copy $B_y$ of $\R$ was inserted into $A$ at the point $(0,y)$ for the Pr\"ufer manifold, $M_\al$ is opened up at $\al$ to make room for $I_\al$ in such a way that for any $\be <\al$ if $\al$ is a limit point of $S_\be$ then every point of $I_\al$ is a limit point of $S_\be$. This will ensure that $M$ is hereditarily separable as we will have ensured that all candidates for countable dense subsets are as dense as needed. At each stage of the induction CH allows use of \cite{B} to deduce that $M_\al$ is homeomorphic to $\B^2$. However $M$ itself cannot be
homeomorphic to $\B^2$; indeed, it is not even metrisable as it is not Lindel\"of.

In his analysis of $\om$--bounded surfaces, Nyikos in \cite[section 5]{N} constructed a tree associated with a significant class of non-metrisable manifolds, those of Type I. A manifold $M$ is of \emph{Type I} provided that there is a sequence $\Sigma =\la U_\al \ra_{\al <\om_1}$ of open subsets of $M$ such that $\bar U_\al$ is Lindel\"of for each $\al$, $\bar U_\al\subset U_{\al +1}$ for each
 $\al$, $U_\lambda =\cup_{\al <\lambda}U_\al$ for each limit $\lambda$, and $M=\cup_{\al <\om_1}U_\al$. Given such a manifold and sequence, the tree of non-metrisable component boundaries associated with $\Sigma$ is the tree $\Upsilon (\Sigma )$ consisting of sets of the form $bd\ C$ where $C$ is a non-metrisable component of $M-\bar U_\al$ for some $\al$ with the order: if $A,B\in\Upsilon (\Sigma)$ then $A\le B$ if and only if $B$ is a subset of a component whose boundary is $A$. Nyikos then uses this tree to enable him to prove the Bagpipe Theorem referred to above.

In \cite{Gr}, Greenwood gives a method of thickening a tree to a manifold so that when we construct Nyikos' tree from the manifold we get the original tree back again (up to isomorphism). By analysing the relationship between the manifold and the tree, she uses this construction to help determine when a non-metrisable manifold contains a copy of $\om_1$ and/or a closed discrete subset. In some Set Theories every non-metrisable Type I manifold contains either a copy of $\om_1$ or an uncountable closed, discrete subspace. On the other hand there are some Set Theories in which there is a non-metrisable Type I manifold which contains neither a copy of $\om_1$ nor an uncountable closed, discrete subspace.

\rk{Acknowledgement}The author is supported in part by a Marsden Fund
Award, UOA611, from the Royal Society of New Zealand.

\Addresses\recd

\end{document}